\documentclass[3p, twocolumn]{elsarticle}

\usepackage{hyperref}
\usepackage{lineno}
\modulolinenumbers[5]

\journal{Journal of \LaTeX\ Templates}

\usepackage{amsmath,amstext,amssymb,graphicx,color,amsthm,float}
\usepackage{subfig, array}

\newcommand{\be}{\begin{equation}}
\newcommand{\ee}{\end{equation}}
\newcommand{\bi}{\begin{itemize}}
\newcommand{\ei}{\end{itemize}}

\usepackage{amsmath}

\DeclareMathOperator*{\argmin}{arg\,min}
    \usepackage{algorithm}
    \usepackage[noend]{algpseudocode}
    \makeatletter
    \def\BState{\State\hskip-\ALG@thistlm}
    \makeatother

\begin{document}

\begin{frontmatter}

 \title{Optimization methods for very accurate Digital Breast Tomosynthesis image reconstruction} 

\author[Disi]{Elena~Morotti\corref{mycorrespondingauthor}}
\cortext[mycorrespondingauthor]{Corresponding author}
\ead{elena.morotti4@unibo.it}

\author[Disi]{Elena~Loli Piccolomini}

\address[Disi]{Department of Computer Science and Engeneering, University of Bologna, Italy}

\begin{abstract}
 
Digital Breast Tomosynthesis  is an X-ray imaging technique that allows a volumetric reconstruction of the breast, from a small number of low-dose two-dimensional projections. 
Although it is already   used in clinical setting, enhancing the quality of the recovered images is still a subject of research. 
Aim of this paper is to propose, in a general optimization framework,  very accurate iterative algorithms for Digital Breast Tomosynthesis image reconstruction, characterized by a convergent behaviour. They are able to detect the cancer object of interest, i.e.  masses and microcalcifications, in the early iterations and to enhance the image quality in a prolonged execution. 
The suggested model-based implementations are specifically aligned to Digital Breast Tomosynthesis clinical requirements and take advantage of a Total Variation regularizer.
We also tune a fully-automatic strategy to set a proper regularization parameter.
We assess our proposals on real  data, acquired from a breast accreditation phantom and a clinical case. The results confirm the effectiveness of the presented solutions in reconstructing breast volumes with particular focus on the masses and microcalcifications.

\end{abstract}


\end{frontmatter}



\section{Introduction \label{intro}}

Digital Breast Tomosynthesis (DBT) is a 3D X-ray cone-beam Computed Tomography (CT)  technique for the early detection of breast tumors \cite{wu2003tomographic,Andersson_2008}. 


While the traditional digital mammography provides a unique 2D breast image, DBT reconstructs the breast as a stack of 2D images by using a comparable radiation dose.
Hence DBT  is also used  in screening programs, because the volumetric reconstruction reduces the tissue overlaps allowing for a better visibility of malignant structures.
DBT is characterized by a limited-angle geometry: since the object is scanned only from a narrow angular range, the DBT projection data is incomplete if compared to classical CT cases. 

The reconstruction algorithm plays an important role, influencing the accuracy of the recovered breast images.
It is well known that traditional fast analytic reconstruction methods, such as  Feldkamp  \cite{Feldkamp_84},  produce poor noisy images in limited-angle tomography, hence they have been left in favour of  Iterative Reconstruction (IR) algorithms \cite{Wu2004a,BLIZNAKOVA201275,sidky_arx2009}. 
IR solvers provide a sequence of solutions, by computing an improved reconstructed volume  at each iteration. 
Many iterative reconstruction solvers have been proposed in literature. An overview of the IR methods is discussed in \ref{sec:stateofart} and a good  paper reviewing IR methods is \cite{Beister2012}.

In this work we consider IR algorithms as solvers of  a model-based formulation  through an unconstrained optimization problem, where the objective function both describes the  CT process  by modelling the physics of the system (including the presence of noise on the projection data) and introduces some image priors. 
Such a mathematical  approach is quite uncommon in 3D  tomographic imaging, where a constrained formulation is preferred \cite{Sidky2009,Choi2010, Ritschl_2011, Luo2017}.
In particular we consider the objective function as the sum of the Least Squares (LS) data fitting term  and the Total Variation (TV) regularization function. The TV regularizer is  chosen by many authors because of its   excellent shape recovering and denoising properties, even if it is known that it can produce staircasing effects when the regularization parameter is too high \cite{LIU2014131,Sidky2009,Choi2010, Ritschl_2011, GraffSidky2016,  Luo2017}.
Hence the choice of the regularization parameter plays a fundamental role in the model-based formulation. 

\begin{figure*}[!t] 	
\centering 
\includegraphics[width=0.9\textwidth]{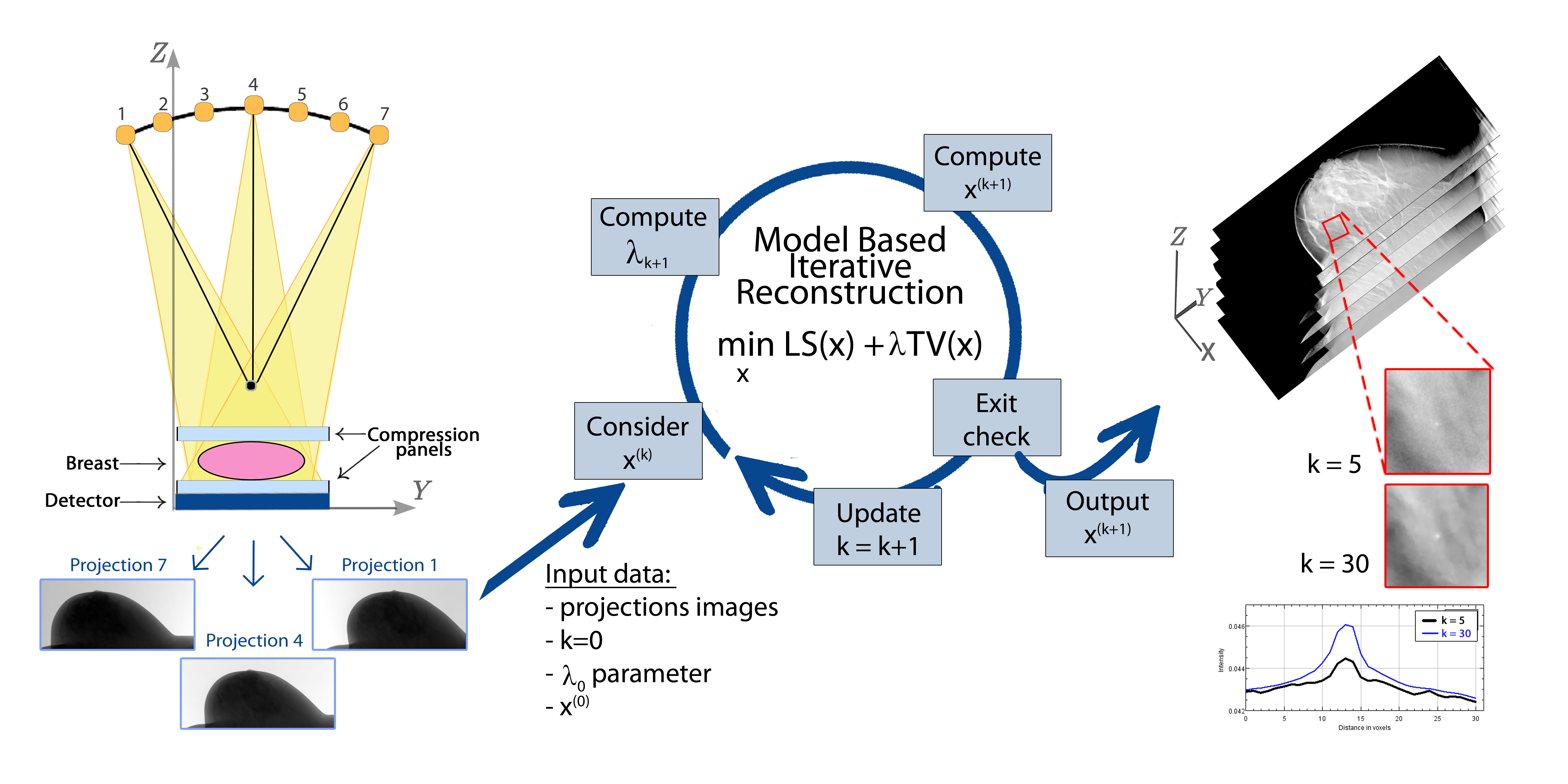} 
\caption{Scheme of the DBT reconstruction process. 
On the left, a draft of the frontal (coronal) section of a  DBT system acquiring  projection images of the breast; in the centre, a chart representing the $k$-th iteration of the algorithm computing the sequence $\{x^{(k)}\}_k$ of approximate solutions by solving the model-based minimization problem; on the right, the evaluation of reconstructed volumes by inspection of cancer objects of interest.}
\label{fig:dbt_sys}
\end{figure*}

Figure \ref{fig:dbt_sys} shows how we approach the entire DBT imaging process, from the numerical modeling of the projection step during the breast scanning, to the reconstructed volume inspection looking for  breast cancer objects, via the implementation of an iterative solver for the model-based minimization problem.

Aim of the paper is to propose both a TV-based  optimization framework and three accurate iterative solvers which use accelerated first order strategies, for DBT image reconstruction.
We are also interested in finding an automatic strategy to set a satisfactory regularization parameter and thus avoid its manually tuning which is infeasible in a clinical setting.

The contribution of this work can be summarized as follows.
\begin{itemize}
    \item We present  three IR solvers in a unique optimization framework which can reconstruct  clinically usable DBT images in few iterations as well as very accurate reconstructions if more iterations are allowed. Even if in clinical routine almost real time reconstructions are required, we remark the importance of improving the image quality with ongoing iterations in longer execution times, for two main reasons:   
    first, having more reliable images can be crucial in difficult diagnosable cases to avoid false responds; 
    second, the fast evolution of multiprocessor boards, such as GPUs, is drastically reducing the time per iteration of the methods, hence we can suppose that more iterations could be performed in clinical reconstructions  in the next future. 
    \item We propose a user independent and computationally effortless rule to set and adapt the regularization parameter at each iteration of the algorithms.
    \item In order to assess our proposals, we implement the methods and test them on real projection data of both a breast accreditation phantom and a human patient.  We analyse the algorithms performance in recovering  the breast tumor objects of interest, by means of measures of merits and visual inspection, at different stages of the iterative reconstruction process. We analyse the volume via its recovered slices, both perpendicularly and  along  the $Z$ direction (see Figure \ref{fig:dbt_sys}).
\end{itemize}
 
The paper is organized as follows.
We present an overview of IR methods in  Section \ref{sec:stateofart}. In Section \ref{sec:model} we state  the optimization framework for the image reconstruction, thus  we illustrate the three proposed IR solvers in Section \ref{sec:methods}. 
Sections \ref{materials} and  \ref{results} present the data sets and the experimental results, respectively. Finally,  Section \ref{concl} contains some conclusions.

 \section{State of art \label{sec:stateofart}}
Iterative approaches have been introduced since the first years of CT, but they have not  been used for long time due to their high computational time request.
Recently, IR methods  got a renewed interest in  scientific communities and among the major vendors,  due to the advent of more performing processors  \cite{Beister2012}. 
As a consequence, a wide amount of IR methods has been proposed to reconstruct  tomographic images and an exhaustive analysis can be found in \cite{GraffSidky2016}.

Initial efforts to solve tomographic imaging with IR methods took an algebraic approach.
Algorithms such as ART, SIRT, SART  and their modifications iteratively solve a linear system of equations by  sequentially projecting a solution onto different hyperplanes \cite{kak1988principles}. 

On the other hand,  the worldwide increasing interest in Compressive Sensing (CS) \cite{Donoho2006} promoted a novel model-based iterative approach, which uses an optimization framework to exploit CS theory.
Among the wide class of model-based IR methods, the so called  Sparsity-Exploiting Image Reconstruction (SEIR) methods have produced significant improvement to the image quality in all the low-dose CT applications (see \cite{LIU2014131, GraffSidky2016} and references therein). 
In particular, many authors introduce  the TV function to take advantage of the sparsity in  the image gradient domain for edge detection \cite{Wu2004a,BLIZNAKOVA201275, Sidky2006a, Sidky2008b, Sidky2008c,Sidky2009b,Lu2010,Chen2013,Ertas2014, Rose2015}.
This property turns into practise as a noise smoothing effect and as a  reliable detection of shape and size of anatomical objects (such as microcalcifications and masses), which are fundamental tasks of DBT imaging.

It is possible to distinguish two main categories of algorithms in the class of SEIR methods: the approximate solvers and the accurate solvers. 
The first one contains algorithms which use, at each step, an algebraic approach (such as SART and SIRT) sequentially and then decrease the TV of the just calculated solution. Examples are  the well-known POCS algorithm and its developments \cite{Sidky2008b, chen2008prior, ramirez2011nonconvex}. They provide reliable reconstructions in few iterations, but the quality of the recovered images strongly depends on the tuning of many inner parameters and the algorithm convergence is not guaranteed.

On the other hand, the accurate solvers are optimization methods which  minimize an objective function defined as a sum of a fit-to-data term and a regularization function. The two quantities are typically weighted by a regularization parameter.
This class is represented by classical optimization methods adapted to the huge size 3D tomographic reconstruction problems. Their solution is proved to converge to the exact solution of the minimization problem.\\
Nowadays, only preliminary investigations on simulations or phantoms  have been  performed to analyse the results of accurate solvers for few-views CT applications \cite{Jensen2012}. 
In our previous works we have investigated a Fixed Point (FP)  algorithm and a Scaled Gradient Projection (SGP) method in \cite{LM16} and \cite{coap2018} respectively, and we have applied them to small simulated data sets. 
A parallel SGP implementation on GPUs has been presented in \cite{Scirep20}, where we were interested in showing the computational  efficiency of the algorithm in terms of execution time.
We also remark that in \cite{park2012fast} an accelerated Gradient Projection method outperformed a sequential reconstruction  algorithm  on data acquired in a sub-sampled 2D  circular geometry.
A further example is the Chambolle-Pock (CP) algorithm which has been applied in  \cite{Sidky2014b} onto 2D circular geometry breast CT, to solve a TV-based convex optimization problem.
In this work we consider the three  iterative optimization  solvers, namely the SGP, FP and CP, to further evaluate their feasibility in reconstructing DBT real volumes and recovering breast tumor objects, like masses and microcalcifications.

Concerning existing rules  for the regularization parameter choice in tomography, in  \cite{siltanen_2016} the authors propose a strategy based on multiresolution and apply it to 2D reconstructions. The proposed rule is very promising, but it  is quite expensive for a very large size 3D application, such as DBT image reconstruction. An  exhaustive list of existing rules for the selection of the regularization parameter  is reported in \cite{siltanen_2016}.

\section{The optimization framework in model-based formulation \label{sec:model}}

Mathematically, tomographic image reconstruction   is an inverse ill-posed problem whose solution can be obtained  by minimizing a suitable objective function related to the physical process.
To define the model describing  image reconstruction  is therefore crucial a deep understanding of the acquisition steps characterizing  the DBT technique. 
A schematic example of a DBT system is shown on the left of Figure \ref{fig:dbt_sys}. 
In DBT routine, the breast is first compressed along the $Z$-axis, over the flat detector plane.
The source moves  along an  arc trajectory  and emits low-dose radiations from a discrete number of angles. Once the X-ray cone-beam has passed  through the body, the detector records its attenuation: the set of the resulting projection images constitutes the raw tomographic data set.
The breast volume to be recovered is composed by a stack of high resolution images, parallel to the detector plane along the $Z$ vertical direction.\\ 
In order to define the numerical model of tomographic image formation, we discretize the 3D object into $N_v$ voxels, whereas the 2D detector panel is made of $N_p$ recording units.
For each fixed projection angle $\theta$ and $i$-th detector recording unit, the Lambert-Beer law relates the projections $P_i^{\theta}$, along a ray $R^{\theta}$, to  the attenuation coefficient function $\mu(w)$  of the voxel $w$ crossed by $R^{\theta}$  \cite{epstein2007introduction} as:
\begin{equation}
  \int_{R^{\theta}} \mu(w)dR=-ln\left(\frac{P_i^{\theta}}{P_0} \right) , \ \ i=1, \ldots, N_p,
  \label{eq:lb}
\end{equation}
where $P_0$ represents the  intensity of the energy emitted by the X-ray source.
The discretization of the integral in \eqref{eq:lb} for  all the  $N_{\theta}$ scanning angles arises the following linear system:
\begin{equation}
    Mx=b.
    \label{eq:sys}
\end{equation}
In equation \eqref{eq:sys} we denote with  $x$ the $N_v$ dimensional vector stacking the attenuation coefficients of all the voxels, while $b$ is the vector of size $N_d=N_p \times  N_{\theta}$ storing all the projections (i.e. the right hand sides of \eqref{eq:lb}) and $M$ is the matrix of size $N_v \times N_d$, built according to the DBT device geometry and representing the projection process onto the detector. 

Some issues arise when solving the linear system \eqref{eq:sys} as an inverse problem, such as  the existence of infinite solutions (since $N_v>N_d$) and the presence of high noise in the reconstructed images (due to the ill-posedness of the problem).
The model-based approach is  introduced  to overcome these numerical controversies, by adding  some a priori information. 
The resulting formulation can be stated as an unconstrained or constrained minimization problem \cite{GraffSidky2016}. 
We consider here the former problem and express it as: 
\begin{equation}
    \min_x f(x) =  J(x)+\lambda R(x)
    \label{eq:minunc}
\end{equation}
where $J(x)$ is a fit-to-data function, $R(x)$ is the prior function (acting here as a regularizer)  and $\lambda$ is the regularization parameter.  

To such DBT mathematical formulation, we can add the box  constraint  $x  \ge 0$ reflecting the non-negativity property of the linear attenuation coefficient $\mu$ in \eqref{eq:lb}.

In particular, in this work we settle $J(x)$ as the Least Squares (LS) function
\begin{equation}
    LS(x)=\|Mx-b\|_2^2
\end{equation}
and $R(x)$ as the  Total Variation (TV) operator defined as \cite{vog02}:
\begin{equation}
    TV(x)= \sum_{i=1}^{N_v} \|\nabla x_i\|_2.
    \label{eq:TV}
\end{equation}
Since TV is not differentiable in the origin, in the algorithms requiring the computation of the gradient, we consider its smoothed version:
\begin{equation}
    TV_{\beta}(x)=\sum_{i=1}^{N_v} \sqrt{\|\nabla x_i\|_2^2+\beta^2}
    \label{eq:tvbeta}
\end{equation}
where  $\beta$ is a small positive parameter \cite{vog02}.
Exploiting the linearity of \eqref{eq:minunc}, the objective function gradient $\nabla f(x) = \nabla LS(x) + \lambda \nabla TV(x)$ can be evaluated by separately computing $\nabla LS(x)$ as
\begin{equation}
    \nabla LS(x) = 2 (M^T M x + M^T b)
\end{equation}
and $\nabla TV(x)$ through finite forward differences.

\section{Iterative optimization methods \label{sec:methods}} 

To solve the  minimization problem \eqref{eq:minunc},  we propose  three accurate solvers  the Scaled Gradient Projection (SGP), the Chambolle-Pock (CP)  and the Fixed Point (FP) methods. For all these methods the convergence to the solution of the model-based  minimization problem d
Among the wide class of optimization methods they have been chosen since they satisfy the requirements necessary to be usable on DBT devices:
 \begin{itemize}
     \item a fast error decreasing in the initial algorithm execution, in order to obtain a good image in few iterations;
     \item a low computational cost per iteration (which is mainly determined by the number of matrix-vector products), to efficiently run the solver in short time;
     \item a limited request of memory, to solve real-size problems on commercially affordable hardware.
 \end{itemize}
A challenging issue (common to the implementation of the three algorithms) is the computation of the projection matrix $M$: since it can not be stored due to its huge dimensions, it must be recalculated at each call.
Thus, this section ends with a focus on the algorithm we use to generate $M$.

 

\subsection{Scaled Gradient Projection  algorithm \label{subsec:SGP}} 

The SGP algorithm is a first order accelerated method. We apply it to solve  the non-negative constrained optimization problem:
\begin{equation}
    \argmin_{x \geq 0} f(x)= LS(x) +\lambda TV_{\beta}(x).
\label{eq:minsgp}
\end{equation}
Algorithm \ref{alg:sgp} reports the main steps of the SGP algorithm.

At each $k$-th iteration, the new solution is computed 
by moving along a descent direction $d^{(k)}$ of a quantity $\eta_k > 0$, as:
\begin{equation}
     {x}^{(k+1)} = {x}^{(k)} + {\eta_k}d^{(k)}.
\end{equation}
The direction $d^{(k)}$ is obtained through a projection ${\cal P}_+$ onto the non-negative orthant:
\begin{equation}
    d^{(k)} = {\cal P}_+\left(\ x^{(k)} -
{\alpha_{k}} S_{k} \nabla f({x}^{(k)}) \right) -  x^{(k)}
\end{equation}
where $\alpha_k$ is the step length and $S_k$ is the scaling matrix (step 7 in Algorithm \ref{alg:sgp}).\\
Essentially, the method follows  a Gradient Projection approach accelerated by choosing the $\alpha_{k}$ step length with Barzilai-Borwein techniques and by introducing a suitable scaling matrix improving the matrix conditioning  \cite{coap2018}.
In particular, the scaling matrix $S_k$ is  a diagonal matrix with entries in a limited interval. 
To update  ${S_{k}}$ (line 5 of Algorithm \ref{alg:sgp}),  we compute a splitting of the objective function gradient into its positive and negative parts, as:
\begin{equation}
\nabla f(x) = V(x) - U(x), 
\label{split}
\end{equation}
where $V(x)>0$ and $U(x) \ge 0 $.
The diagonal elements ${s}^{(k)}_{j,j}$ of $S_{k}$ 
are updated, for $j=1, \ldots N_v$ as:
\begin{equation}
s^{(k)}_{j,j} =  \min\left(\! \rho_{k}, \max\left(\! \frac{1}{\rho_{k}}, \frac{x^{(k)}_j}{V_j(x^{(k)})} \right)\right)
\end{equation}

where $\{ \rho_{k} \}_k$ is a decreasing positive sequence.

Regarding the convergence, it is proved in  \cite{BP15} that the SGP algorithm converges without any further restriction on the step length $\alpha_k$ and on the scaling matrix $S_k$ to the unique minimum of \eqref{eq:minsgp}.
In \cite{BPR16}, the authors proved that the theoretical convergence rate  of the SGP method is $\cal{O}$(1/k).


\begin{algorithm}
\caption{Scaled Gradient Projection algorithm (SGP)}\label{alg:sgp}
\begin{algorithmic}[1]
\Require $M, b, \lambda $
\State {\bf Initialize: } ${x}^{(0)}\ge 0, \ \ \gamma,\sigma\in(0,1),\ \ 0 < \alpha_{min}\le\alpha_{max},$
\State k=0 
\While {not convergence}    
    \State{Compute $g^{(k)} = 2(M^TMx^{(k)} + M^Tb) + \lambda \nabla TV_\beta (x^{(k)})$ }
    \State{Compute  ${S_{k}}\in {S_{\rho_{k}}}$}
    \State{Define  ${\alpha_{k}}\in [\alpha_{min}, \alpha_{max}]$ with alternate BB rules}
    \State $ d^{(k)} = {\cal P}_+\left( x^{(k)} -
    		{\alpha_{k}} S_{k}
    		g^{(k)} \right) -  x^{(k)}$
    \State $\eta_{k} = 1$
    \While {$ 		f({x}^{(k)} + {\eta_{k}} d^{(k)})  > f({x}^{(k)}) +
    		\sigma{\eta_{k}} (g^{(k)})^T{d}^{(k)}	$ }
        \State ${\eta}_{k} = \gamma {\eta}_{k}$
    \EndWhile
    \State ${x}^{(k+1)} = {x}^{(k)} + {\eta}_{k} d^{(k)}$
    \State k = k+1
\EndWhile    
\Ensure $x^{(k)}$
\end{algorithmic}
\end{algorithm}

\subsection{The Fixed Point algorithm  \label{subsec:FP}} 

The FP algorithm for the solution of the minimization problem:
\begin{equation}
    \argmin_x f(x)= LS(x) +\lambda TV_{\beta}(x)
    \label{eq:minu}
\end{equation}
has been firstly proposed for image denoising by Rudin, Osher and Fatemi in \cite{RUDIN1992}.
 Starting from this approach we derived the lagged diffusivity FP Algorithm \ref{alg:fp} for 3D tomographic image reconstruction. \\
\begin{algorithm}
\caption{Lagged diffusivity Fixed Point algorithm (FP)}\label{alg:fp}
\begin{algorithmic}[1]
\Require $M, b, \lambda, maxiter$
\State {\bf Initialize: }${x}^{(0)} \ge 0$ 
\For{$k=0$ to $maxiter-1$}    
    \State{Compute $g^{(k)} = 2(M^TMx^{(k)} + M^Tb) + \lambda \nabla TV_\beta (x^{(k)})$ }
    \State Solve the linear system $H_k d^{(k)} = -g^{(k)} $, where $H_k=M^TM+\lambda L(x^{(k)})$, with the Conjugate Gradient method.
    \State ${x}^{(k+1)} = {x}^{(k)} + d^{(k)} $        
\EndFor     
\Ensure $ {\cal P}_+(x^{(k+1)})$
\end{algorithmic}
\end{algorithm}

At each $k$-th iteration, the FP algorithm updates the  solution with the following rule:
\begin{equation}
    {x}^{(k+1)} = {x}^{(k)} + d^{(k)} 
\end{equation}
where  the descent direction $d^{(k)}$ is computed by solving a linear system $H_k d^{(k)}=-\nabla f({x}^{(k)})$  (line 4 of Algorithm \ref{alg:fp}).
The matrix  $H_k=M^TM+\lambda L(x^{(k)})$ in line 4 approximates the Hessian matrix. It contains in fact the seven diagonals banded matrix $L(x^{(k)})$ which is the discretization matrix of the diffusion operator $L(x)$ so that  $L(x)x = \nabla TV(x)$ \cite{vog02}.
We solve the linear system with very few iterations of a Conjugate Gradient (CG) algorithm \cite{hestenes}: we stop it far before convergence, both to limit the computational time and to prevent  noise from affecting the solution.
We remark that each CG iteration requires a matrix-vector product involving $H_k$ and that, to save memory space, we perform it without storing the matrix $H_k$: we only store $L(x^{(k)})$ and re-compute $M$ and $M^T$ at run time. 
At the end, we project the last computed solution onto the non-negative orthant.
For more details on the FP method applied to tomographic image reconstruction and its convergence, see \cite{LM16} and \cite{chan1999convergence} respectively.

\subsection{The Chambolle-Pock  algorithm  \label{subsec:CP}} %


The CP algorithm has been firstly proposed in \cite{ChambolleP11}  for the solution of the general minimization problem:
\begin{equation}
\argmin_x f(x) =  F(Kx)+G(x),
\label{eq:cp}
\end{equation}
where $G$ is a convex, lower-semicontinuous and proper function, $F$ is convex and lower-semicontinuous and K is a continuous linear operator. 
To fit the problem statement \eqref{eq:cp} and fully exploit the linearity of $K$, we assign:
\begin{equation}
    F(Kx) = LS(x) + { \lambda} TV(x)
\end{equation}
which results in defining the operator $K$ as a matrix composed by the four following blocks:
\begin{equation}
    K = \begin{pmatrix} M \\ \nabla_x \\ \nabla_y \\ \nabla_z. \\ \end{pmatrix} 
\end{equation}
where $\nabla_x,  \nabla_y$ and $\nabla_z$ are the forward differences operators acting along the $X, Y$ and $Z$ axes respectively. 
In order to include the non-negative constraints, we fixed $G(x)=\delta_{\Omega}(x)$ as the indicator function of the convex set $\Omega=\{ x: x \geq 0 \}$, i.e. 
\begin{equation}
    \delta_{\Omega}(x)=\begin{cases}
    0 & x \in \Omega \\
    \infty & x \notin \Omega. \\
    \end{cases}
\end{equation}

Considering the convex conjugate $F^*$ of $F$, defined as $F^*(y) = \max_x \{ x^Ty - F(x) \} $, and the proximal mappings of $G$ and $F^*$, i.e.
\begin{equation}
\begin{split}
prox_\sigma [F^*](y) & =  \argmin_{\bar y} \Big\{ F^*(\bar y) + {1 \over 2\sigma} \|y-\bar y\|_2^2  \Big\}      \\
prox_\tau [G](x)     & =    \argmin_{\bar x} \Big\{G(\bar x) + {1 \over 2\tau} \|x-\bar x\|_2^2  \Big\}  ,     
\end{split}
\label{eq:prox}
\end{equation}
the $k$-th CP iteration can be described with the following three steps:
\begin{enumerate}  
    \item compute $y^{(k+1)}$ as $prox_\sigma [F^*](y^{(k)} + \sigma K \bar x^{(k)}) $;
    \item compute $x^{(k+1)}$ as $prox_\tau [G] (x^{(k)} - \tau K^T y^{(k+1)}) $;
    \item define $\bar x^{(k+1)}$ with an extrapolation step:
  $\bar x^{(k+1)}=x^{(k+1)} + \theta ( x^{(k+1)} - x^{(k)}) $ and $\theta>0$.
\end{enumerate}

In particular, the proximal mapping $prox_\sigma [F^*]$ can be computed as sum the two independent blocks, as in lines 5-6 and 7-8 of the Algorithm \ref{alg:cp}; its detailed derivation  can be found in \cite{Sidky2014b}. 
The proximal mapping of $G$ is defined as:
\begin{equation}
\begin{split}
prox_\tau [G](x) & =    \argmin_{\bar x} \Big\{ \delta_{\Omega}(x) + {1 \over 2\tau} \|x-\bar x\|_2^2  \Big\}  \\
     & = \argmin_{\bar x \in \Omega} \Big\{  {1 \over 2\tau} \|x-\bar x\|_2^2  \Big\} \\
     & = {\cal P}_{+}(x)
\end{split}
 \label{eq:cpProxG}
\end{equation}
hence it is exactly the projection $ {\cal P}_{+}(x)$ of $x$ onto the feasible set $\Omega$ (lines 9-10 of Algorithm  \ref{alg:cp}).
The  updated iterate $x^{(k+1)}$ is computed with a FISTA strategy as in line 11 of  Algorithm \ref{alg:cp}.
The algorithm convergence is demonstrated in \cite{ChambolleP11}.\\
We finally observe that the algorithm needs to compute the value $ \Gamma$ (line 1 of  Algorithm \ref{alg:cp}):  to estimate the matrix 2-norm as $ \Gamma \approx \|K\|_2 = \sqrt{\rho(K^T K) }$ (where $\rho$ is the spectral radius of a matrix), we perform two iterations of the power method for the maximum eigenvalue computation \cite{golub2012matrix}.   

\begin{algorithm}
\caption{Chambolle Pock algorithm (CP)}\label{alg:cp}
\begin{algorithmic}[1]
\Require $M, b, \epsilon, maxiter$
\State {\bf Compute: }  $\Gamma$ as an approximation of  $\|K\|_2$
\State {\bf Initialize: } $\tau=\sigma=\frac{1}{\Gamma}>0$, $\theta \in [0,1]$
\State {\bf Initialize: } ${x}^{(0)} \ge 0, \bar x^{(0)}, {y}^{(0)}$ and ${w}^{(0)}$  to zeros-vectors
\For{$k=0$ to $maxiter-1$}    
    \State $\bar y^{(k)}=y^{(k)}+ \sigma(M \bar x^{(k)}-b)$
    \State $y^{(k+1)} = max(\|\bar y^{(k)}\|_2 -\sigma \epsilon) {\bar y^{(k)} \over \|\bar y^{(k)}\|_2}$
    \State $ \bar w^{(k)}=w^{(k)}+ \sigma (\nabla_x,\nabla_y,\nabla_z) \bar x^{(k)}$
    \State $ w^{(k+1)}=\bar w^{(k)}(\lambda /max(\lambda, |\bar w^{(k)}|)$ 
    \State $ x^{(k+1)} =  x^{(k)}-\tau(M^T y^{(k+1)}+(\nabla_x,\nabla_y,\nabla_z)^T w^{(k+1)}$ 
    \State $ x^{(k+1)}={\cal P}_+( x^{(k+1)})$ 
    \State $\bar x^{(k+1)}=x^{(k+1)} + \theta (x^{(k+1)}-x^{(k)})$
\EndFor     
\Ensure $x^{(k+1)}$
\end{algorithmic}
\end{algorithm}

\subsection{User-independent choice of the regularization parameter \label{subs:TeoriaParamreg}  }

In model-based optimization approach \eqref{eq:minunc}, the choice of the regularization parameter $\lambda$  plays a key role for the quality of the reconstruction and it represents a crucial challenge in a clinical setting, where the trial-and-error approach is not doable for each reconstruction.
Moreover, experimental results show that 
a strong regularization is required in the first iterations, to avoid noise propagation and force the algorithm towards a good solution, whereas a  weaker regularization in the last iterations can prevent the TV staircaising effects on final reconstructions.
Hence, we propose to reduce the regularization weight along the iterations, by choosing the $\lambda$ values with a decreasing updating rule.
Interestingly, state of art studies have already proposed semi-automatic rules for the selection of a decreasing sequence $\{ \lambda_k \}_k, \ k = 1, \dots$ of regularization parameters defining a sequence of minimization problems stated as \eqref{eq:minu}, whose solutions converge to a good reconstructed image. 
See \cite{noi_ip2019} for more details and the convergence proof.

We propose the following fully-automatic strategy to compute a decreasing sequence $\{ \lambda_k \}_k$.
At the beginning of our algorithm, we leave out the regularization by setting the first parameter $\lambda_0=0$: we are in fact interested in a very good data fitting, to recover as many image features as possible.
Next, the starting value $\lambda_1 $ is set to balance the residual norm and the amount of TV of the first iterate. 
Afterward, we propose to  decrease $\lambda$ of a constant factor $1/k$ at each $k$-th  iteration, since we need a very simple and computationally cheap rule, reducing the regularization weight slightly. 
The resulting strategy is summarized in the following scheme and it can be introduced  in each of the previously considered algorithms.

\begin{itemize}
    \item Set $\lambda_0 = 0$ to initialize the algorithm and run the first iteration (labelled with k=0) to compute $x^{(1)}$;
    \item Set $ \lambda_1 = { \sqrt{LS(x^{(1)}) }  \over TV (x^{(1)}) } $ and use it to compute $x^{(2)}$; 
    \item For each $k \ge 2$, set
    \begin{equation} \label{eq:UpdateLambda}
        \lambda_{k} = {1 \over k}  \lambda_1
    \end{equation} and use it to compute  $x^{(k+1)}$.
\end{itemize}


\subsection{The projection matrix algorithm  \label{subsec:DD}} 
Besides the choice of the model parameter and the solver, in optimization approach a key point consists in  numerical modeling  the geometric projection process, schematically  displayed from a frontal view in Figure \ref{fig:dbt_sys}, through a matrix.  

\begin{figure}[!t] 	
\centering
\subfloat[]{\includegraphics
[width=0.22\textwidth]{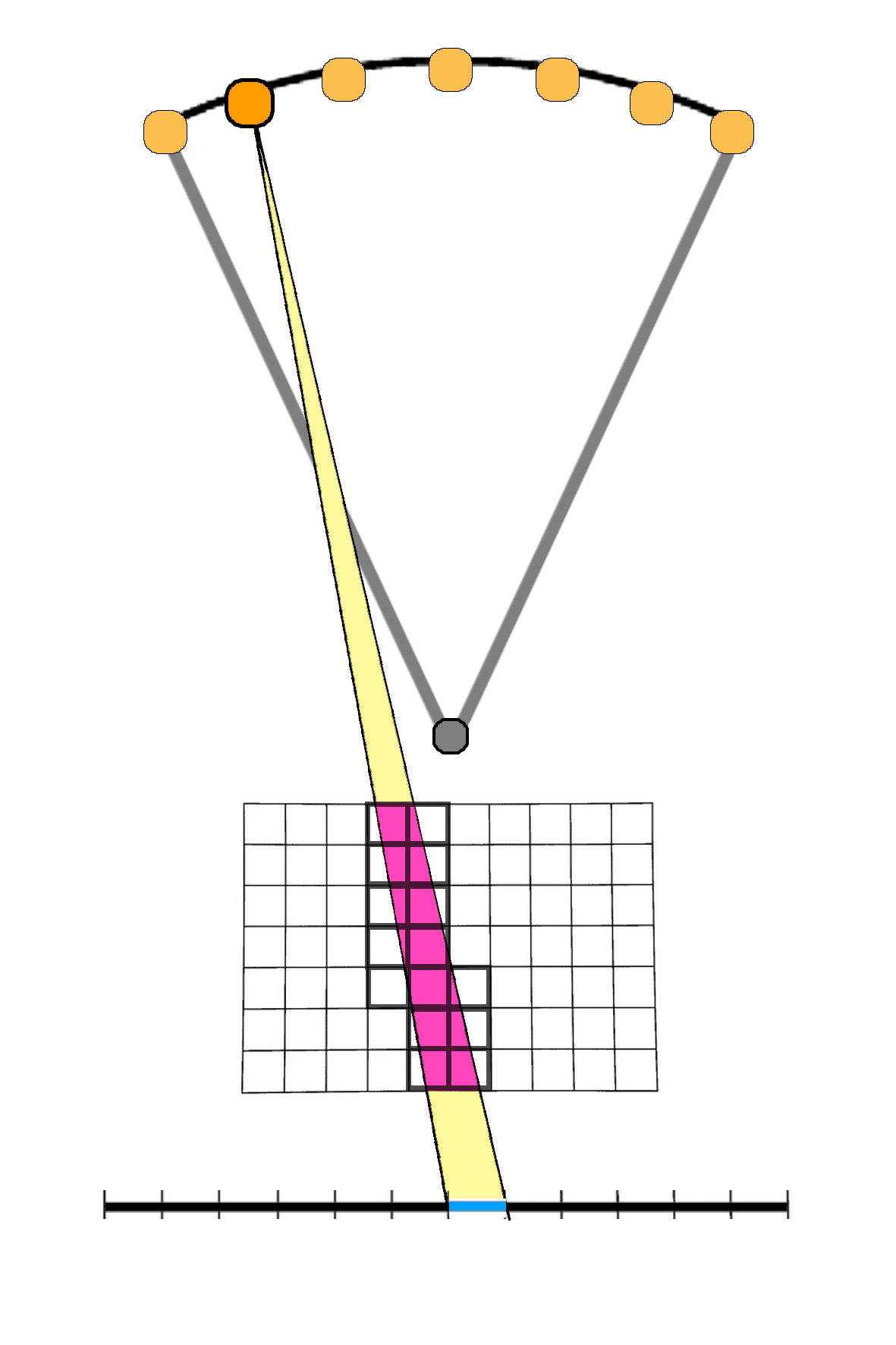} \label{fig:DD_1}}
\subfloat[]{\includegraphics[width=0.2\textwidth]{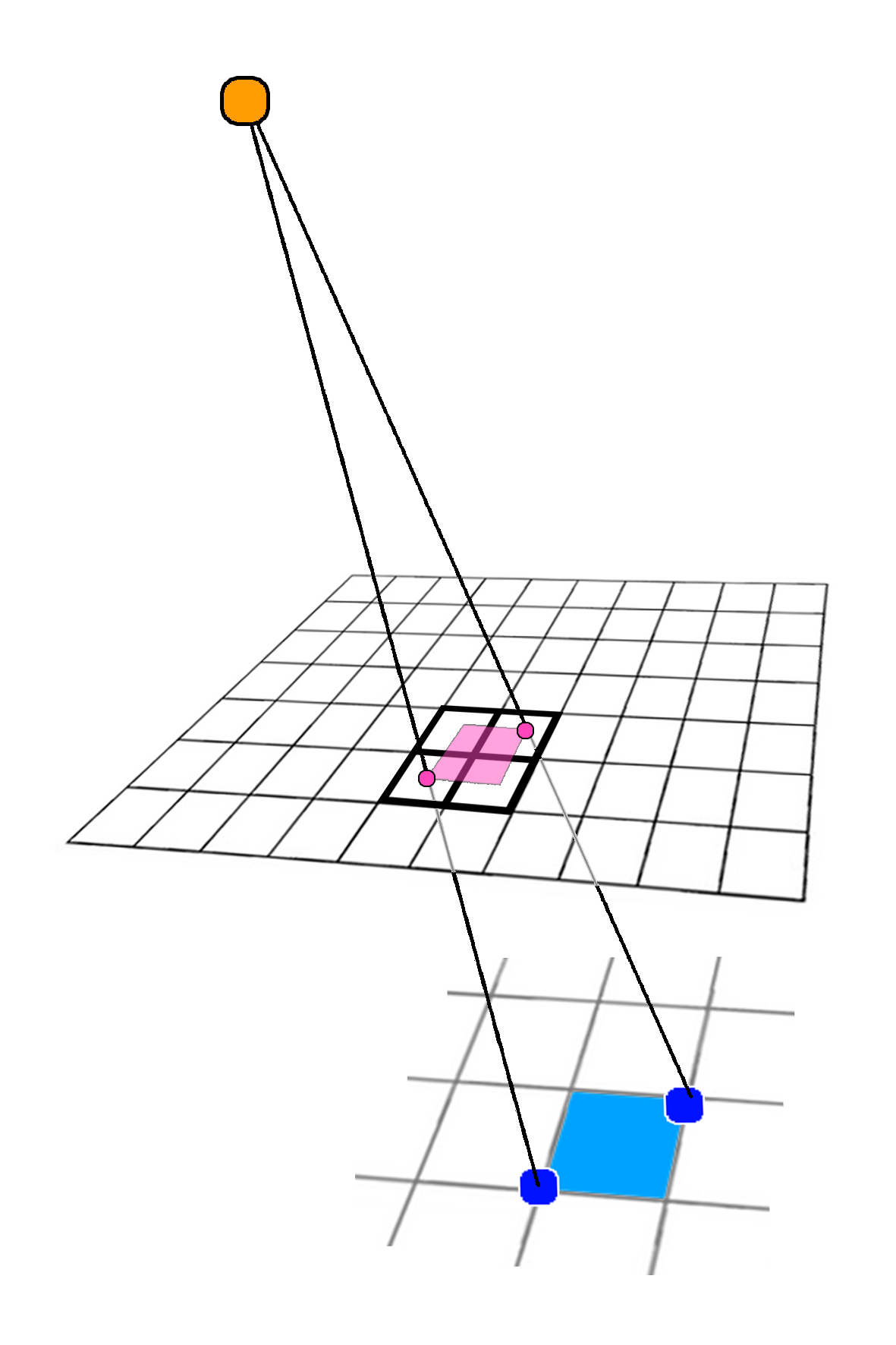} \label{fig:DD_2}}
\hfil\\
\caption{Schematic draw representing the Distance Driven approach to compute the system matrix.
(a) View on the $YZ$ plane of an X-ray projection onto a single pixel, from a fixed angle. The intersection of the X-ray beam with the volume is highlighted in magenta.
(b) The magenta area represents the backward projection of the blue recording unit onto a volume slice  parallel to the $XY$ plane.  }
\label{fig:DD}
\end{figure}

The coefficient  matrix $M$ of the linear system  \eqref{eq:sys} is commonly called \emph{projection operator} in tomography, since it represents the action of the tomographic system in projecting an object onto the detector,
whereas the matrix modeling the  backprojection of  the tomographic data onto a volume is called \emph{backprojection operator}.  In the proposed optimization algorithms, the backprojection coincides with the transpose matrix $M^T$.\\
Different algorithms have been proposed in literature for the computation of the matrix $M$.  We have adopted the Distance Driven (DD), which  accurately models the discretization of the Lambert-Beer's law \eqref{eq:lb} for cone-beam projections \cite{DeManBasu2004}. 
In DD, $M$ of size $N_v \times (N_p \times N_{\theta})$ is constituted by $N_\theta$ submatrices $M^\theta$ of size $N_v \times N_p$. Each element $M^\theta_{i,j}$  represents the contribution of the $j$-th voxel (for $j=1, \ldots N_v$) to the projection  onto the $i$-th detector pixel (for $i=1, \ldots N_p$), for a projection angle $\theta$.
Images in Figure  \ref{fig:DD} help in understanding the DD procedure.
In Figure \ref{fig:DD} (a), for a scanning angle, we consider the X-ray cone-beam projecting onto the $i$-th blue pixel and intersecting  the voxels with bold contours (in the magenta coloured area). Only these voxels contribute to the value of the projection in the considered pixel.
In Figure \ref{fig:DD} (b) we highlight the i-th cell of the detector (the blue area)  and its backward footprint on a plane parallel to the detector (the magenta area).  The ratio between the magenta area inside the  $j$-th voxel and the whole magenta extension  is proportional to the value $M_{i,j}$ of the matrix. 
For all the  voxels $j$ not contributing to the $i$-th projection the corresponding matrix element  $M_{i,j}=0$; hence $M$ is extremely sparse.
However, despite the  huge number of nonzero elements, for its very large size, in real applications $M$  cannot be stored and it must be recomputed  whenever a matrix-vector product is needed.

We finally remark that we have  modified the general approach presented in \cite{DeManBasu2004} by efficiently  exploiting the characteristics of our specific mammographic setting.
Really,  since the DBT detector is a stationary flat panel and it is parallel to the compression plane of the breast, the footprints   can be directly projected onto the detector plane, thus avoiding  the use of an intermediate projection plane and  further computational costs.

\section{Materials \label{materials}}

\subsection{DBT system configuration  \label{subs:Geom}}

Our tests are performed on  the digital system {\it Giotto Class} of the Italian I.M.S. Giotto Spa company in Bologna \cite{siteIMS}. 
The source executes $N_{\theta}=11$  scans from equally spaced angles in  approximately $30$ degrees  range; in the highest vertical position, the source is about $70 cm$ over the detector.
The stationary digital detector has a sensitive area of $24 \ cm \times  30 \ cm$  and squared pixel pitch of 0.085 mm; 
the reconstructed voxel dimensions along the three cartesian axes are  $\Delta_x = \Delta_y = 0.090\ mm$ and $\Delta_z = 1\ mm$ respectively.\\
The system uses a polychromatic ray with energies in a narrow range around 20 keV to avoid the photon scattering. As always happens in CT reconstruction algorithms, we approximate the polychromatic beam with a monochromatic one.

\subsection{Data sets \label{subs:Phantom}}

We consider two data sets in our experiments: a breast 3D phantom and a clinical acquisition from a human subject.
Both volumes contain  the objects of interest for breast cancer detection, i.e. small high contrast microcalcifications and larger but lower contrasted masses.

The  phantom is the model 020 of BR3D breast imaging phantom, produced by CIRS Tissue Simulation and Phantom company \cite{SitoMortadella}. 
It is characterized by a heterogeneous background, where adipose-like and gland-like tissues are mixed in about 50/50 ratio and it is made of six slabs that may be arranged to create multiple anatomical backgrounds.
Each slab has a semicircular shape and its size is  $10\ cm \times 18 \ cm$.
Inside one of  them, we find acrylic spheres simulating breast masses (MSs), 1 $cm$ length fibers and many clusters of calcium carbonate specks simulating microcalcifications (MCs). 
We report in Table \ref{tab:SpecificheMortadella} the length of the diameters of all the MSs and of each sphere of a  MC cluster. 
In particular, we reconstruct a volume of 50 slices of $11.4\ cm \times 21 \ cm$ and we analyze in the reconstructed images objects with different diameters, such as the microcalcifications in clusters 3, 5 and  6 (having  230, 165 and 130 $\mu m$ diameter, respectively) and the second and fourth mass (with diameter 4.7 and 3.1 $mm$,  respectively). 
All such objects lye on the same slice.

As human DBT data set, we have chosen a case containing  microcalcifications, circular and spiculated masses. The clinical volume  is constituted of 55 slices of $10.5 \ cm \times 20.7 \ cm$.


\begin{table}[h]
    \centering
    \begin{tabular}{c|c|c|c|c|c|c}
      & 1& 2& 3 & 4 & 5& 6\\
      \hline
       MC  & 400 & 290 & 230 & 196 & 165 & 130 \\
       MS  & 6300 & 4700 &3900 & 3100 & 2300 & 1800\\ 
    \end{tabular}
    \caption{Diameters of a microcalcification (MC) in a cluster  and of a mass (MS) in the BR3D phantom as reported in \cite{SitoMortadella}. Measures are in micrometers ($\mu m$). }
    \label{tab:SpecificheMortadella}
\end{table}

\subsection{Measure and graphics of merits}
In order to quantitatively evaluate the reconstructed objects of interest in the volumes, we compute two widely used  measure of merits: the Contrast-to-Noise Ratio (CNR) and the Full-Width at Half Maximum (FWHM).

The CNR measure on a mass is calculated as:
\begin{equation}
CNR_{MS} = \frac{\mu_{MS} -\mu_{BG}}{\sigma_{MS} -\sigma_{BG}}
\label{eq:cnrms}
\end{equation}
where $\mu$ and $\sigma$ are the mean and standard deviation computed on the reconstructed volume, in small regions located inside the mass (MS) or in the background (BG).
Similarly, we define the CNR measure on a microcalcification as:
\begin{equation}
CNR_{MC} = \frac{\mathcal{M}_{MC} -\mu_{BG}}{\sigma_{BG}}
\label{eq:cnrmc}
\end{equation}
where $\mathcal{M}_{MC}$ is the maximum intensity  inside the considered microcalcification (MC).
Higher values of the CNR indices reflect a better detection of an object from the background.

To compute the FWHM parameter, we consider the transverse slice (parallel to the $XY$ plane) where the microcalcification lies and we extract the  \emph{Plane  Profile (PP)} along the $Y$ axes. The FWMH index is computed as:
\begin{equation}
    FWHM=2 \sqrt{2 \ln{(2)}} d
    \label{eq:fwhm}
\end{equation}
where $d$ is the standard deviation of the gaussian curve fitting the PP. 
We remark that 
\begin{equation}
w=FWHM \cdot  \Delta_y
\label{eq:w}
\end{equation}
approximates the width of the examined microcalcification.
The Plane Profiles are also useful tools to evaluate the reconstruction accuracy on the transverse plane.

To estimate the solver effectiveness along the  $Z$ direction, which is the most challenging purpose in DBT imaging, we plot the Artifact Spread Function (ASF) vector, whose components are computed on a microcalcification as: 
\begin{equation}
 ASF(z)=\frac{|\mu_{MC}(z)-\mu_{BG}(z)|}{|\mu_{MC}( \bar z)-\mu_{BG}(\bar z)|}, \ \ \forall z=1, \ldots , N_z
\label{eq:ASF}
\end{equation}
where $\mu(z)$ is the mean of the reconstructed values inside a circular region of three pixels diameter inside the considered MC and in the background,  $\bar z$ corresponds to the slice where the object is on focus and $N_z$ is the total number of discrete slices. Similarly, we compute the ASF for the masses.

\section{Numerical results and discussion\label{results}   }

In this section we present the results obtained with the proposed optimization approach and the accurate solvers described in section \ref{sec:methods}. 
At first we compare the BR3D phantom reconstructions produced by the  SGP, FP and CP solvers in a similar computational time.
Then, to analyse the best obtainable image  quality we have run the SGP up to convergence both on the phantom and on the clinical data set.
In all the  previous tests we used  a constant value of $\lambda$, set by trial and error. 
At last,  we test the automatic rule proposed in section \ref{subs:TeoriaParamreg} to decrease the $\lambda$ values along the iterations.

\subsection{Methods comparison for early reconstructions}\label{subs:comparison}
Aim of this paragraph is to show the behaviour of the proposed solvers at different stages of their executions.
We fixed 5 and 15 iterations: the workload of 5 iterations is compatible with the execution of a reconstruction on a commercial hardware  in a clinical setting, while in 15 iterations we get  fairly accurate  reconstructions with all the three methods, reflecting that they are sufficiently close to the convergence solution. 
Each SGP and CP iteration requires approximately the same time, whereas in the special case of FP solver, the number of allowed iterations  corresponds to the sum of the external  and  CG  iterations. Since we perform 4 CG iterations, we have stopped  the  FP algorithm  after one or three outer iterations (loop $k$ in Algorithm \ref{alg:fp}), respectively.

The value of $\beta$ in \eqref{eq:tvbeta} has been fixed as $\beta=0.001$.
Since  the three considered algorithms  solve a slightly different optimization problem, the three $\lambda$ parameters have been chosen  independently for each method to achieve the best reconstruction in 5 iterations: we have set $\lambda=0.005$ for both SGP and CP methods and $\lambda=0.001$ for the FP algorithm.

\begin{figure}[!t] 	
\centering 
\subfloat[SGP]{\includegraphics[width=0.23\textwidth]{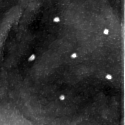} }
\subfloat[SGP]{\includegraphics[width=0.23\textwidth]{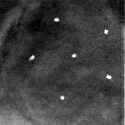}}
\hfil \\
\subfloat[FP]{\includegraphics[width=0.23\textwidth]{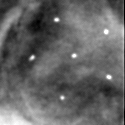} }
\subfloat[FP]{\includegraphics[width=0.23\textwidth]{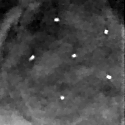} }
\hfil \\
\subfloat[CP]{\includegraphics[width=0.23\textwidth]{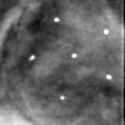} }
\subfloat[CP]{\includegraphics[width=0.23\textwidth]{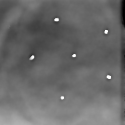} } 
\hfil
\caption{Reconstructions of microcalcification cluster number 3 in BR3D phantom obtained with SGP, FP and CP methods. 
On the left column, reconstructions in 5 iterations; on the right, reconstructions in 15 iterations.}
\label{fig:confronto3metodi}
\end{figure}

\begin{figure*}[!t] 	
\centering 
\subfloat[PP in 5 iterations] {\includegraphics[width=0.47\textwidth]{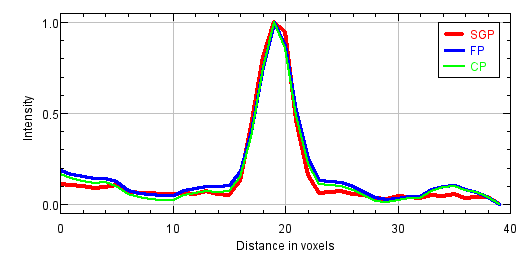} \label{fig:Vp4}}
\subfloat[ASF  in 5 iterations] {\includegraphics[width=0.47\textwidth]{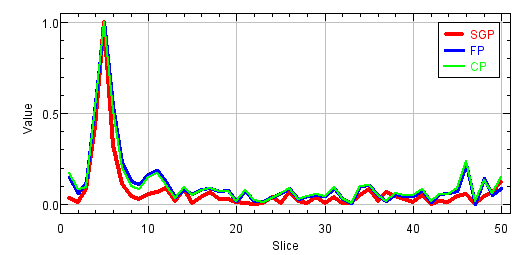}  \label{fig:Dp4}}
\hfil \\
\subfloat[PP in 15 iterations] {\includegraphics[width=0.47\textwidth]{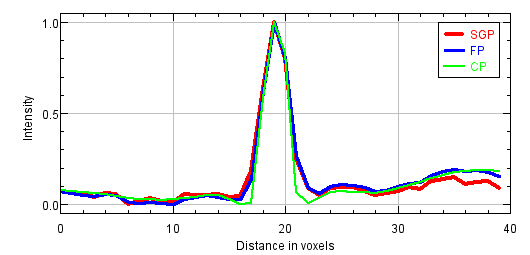} \label{fig:Vp12}}
\subfloat[ASF  in 15 iterations] {\includegraphics[width=0.47\textwidth]{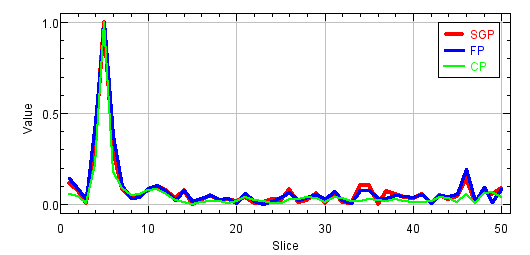} \label{fig:Dp12}}
\hfil
\caption{Plots of the Plane Profile on the left and of the ASF vectors on the right, taken over one microcalcification of cluster number 3 in BR3D phantom obtained.
In all the plots: the red line corresponds to  SGP method, the  blue line to FP method and the green line to CP method.
}
\label{fig:curveConfronto3metodi}
\end{figure*}

In the following analysis, we focus on  the reconstruction of MC cluster number 3 in  the BR3D phantom. 
For each solver, Figure \ref{fig:confronto3metodi} reports a $125 \times 125$ pixels crop taken from the fifth slice of the reconstructions in 5 and 15 iterations. 
The images are represented by automatically enhancing the gray level contrast computed on the same considered region.
In Figure \ref{fig:curveConfronto3metodi} we  compare the PP and ASF curves taken on one MC.

Looking at  Figure \ref{fig:confronto3metodi},  we observe that the detection of the MC cluster is comparable at equal iterations whereas the background appears slightly different for the three methods. For example,  in the case of CP reconstruction in 15 iterations it looks smoother and more blurred. 
Focusing on the objects of interest, we notice that in 5 iterations the MCs are perfectly visible;
moreover, in 15 iterations the MC edges are sharper as confirmed by the plots (a) and (c)  in Figure  \ref{fig:curveConfronto3metodi}.
From plots (b) and (d) of Figure \ref{fig:curveConfronto3metodi} we observe that in all the three reconstructions the object is  placed in the correct slice and it is not diffused in the adjacent layers.    
Hence we can conclude that the proposed model-based optimization framework yields good quality images in early reconstructions, regardless the applied solver.

\subsection{ SGP algorithm insights \label{subs:SGP_phantom}}

\begin{figure}[!t] 	
\centering
\includegraphics[width=0.45\textwidth]{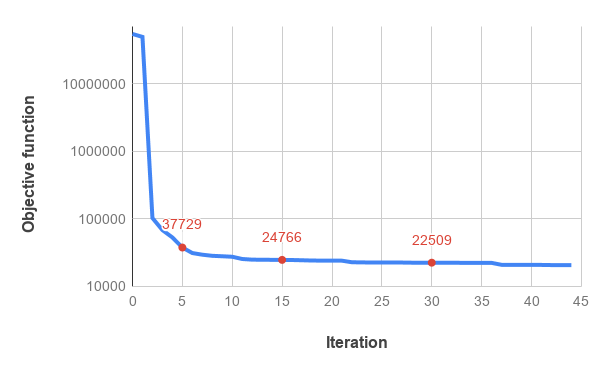} 
\caption{Objective function values vs. iteration number for the SGP execution on the phantom test. The convergence has been reached after 44 iterations by satisfying condition \eqref{eq:stop}. The red labels outline the function values at 5, 15 and 30 iterations. }
\label{fig:ObjFuncMortadella}
\end{figure}

In the following, we raise the SGP as the  representative solver in the proposed optimization framework. 
We  explore the performance of the optimization approach on many different objects  of the BR3D phantom (such as microcalcifications of very small diameter and masses with a low contrast with the background tissue) and we analyse the quality of the reconstructions also after 15 iterations, i.e. approaching convergence.

\begin{figure}[!t] 	
\centering
\subfloat[5 iterations]{\includegraphics[height=0.21\textwidth, trim={42.8cm 6cm 34.7cm 34cm}, clip]{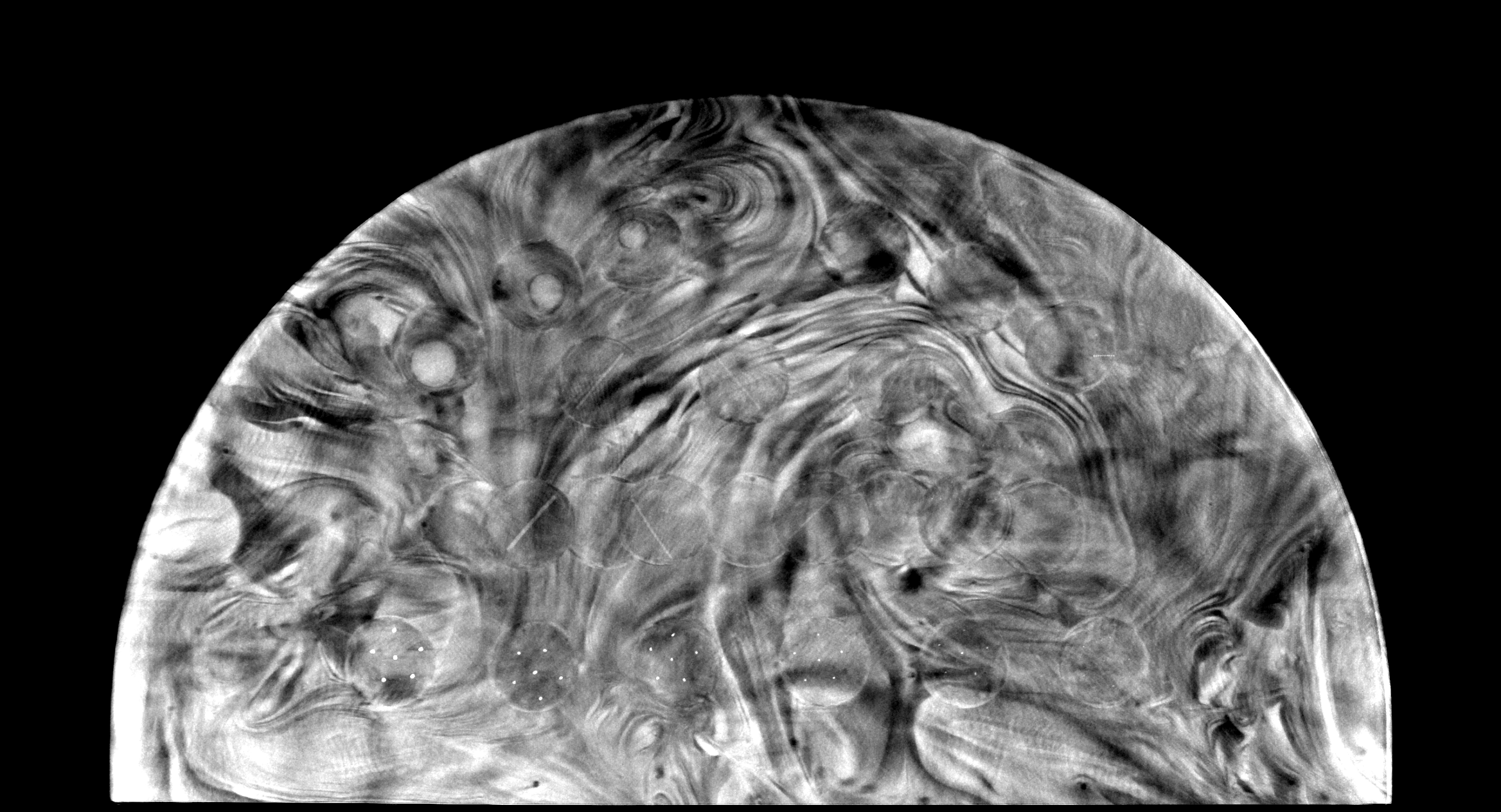} \label{fig:MicroMortadella5}}
\subfloat[5 iterations]{\includegraphics[height=0.21\textwidth, trim={27cm 26.3cm 50.5cm 13.5cm}, clip]{Sgp_RECO_4_Layer7.png} \label{fig:MassaMortadella5}} 
\hfil\\
\subfloat[15 iterations]{\includegraphics[height=0.21\textwidth, trim={42.8cm 6cm 34.7cm 34cm}, clip]{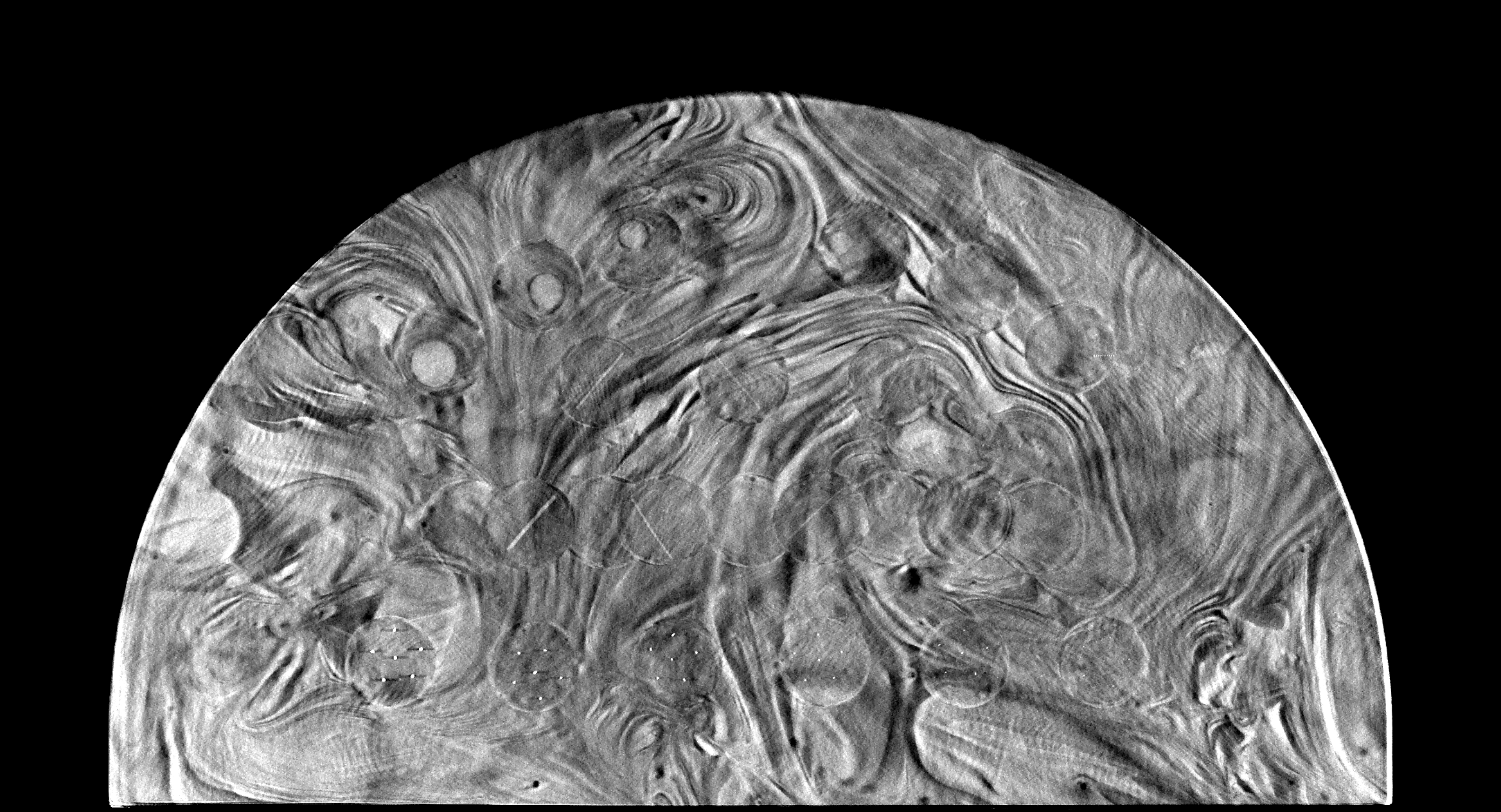} \label{fig:MicroMortadella15}}
\subfloat[15 iterations]{\includegraphics[height=0.21\textwidth, trim={27cm 26.3cm 50.5cm 13.5cm}, clip]{Sgp_RECO_14_Layer7.png} \label{fig:MassaMortadella15}}
\hfil\\
\subfloat[30 iterations]{\includegraphics[height=0.21\textwidth, trim={42.8cm 6cm 34.7cm 34cm}, clip]{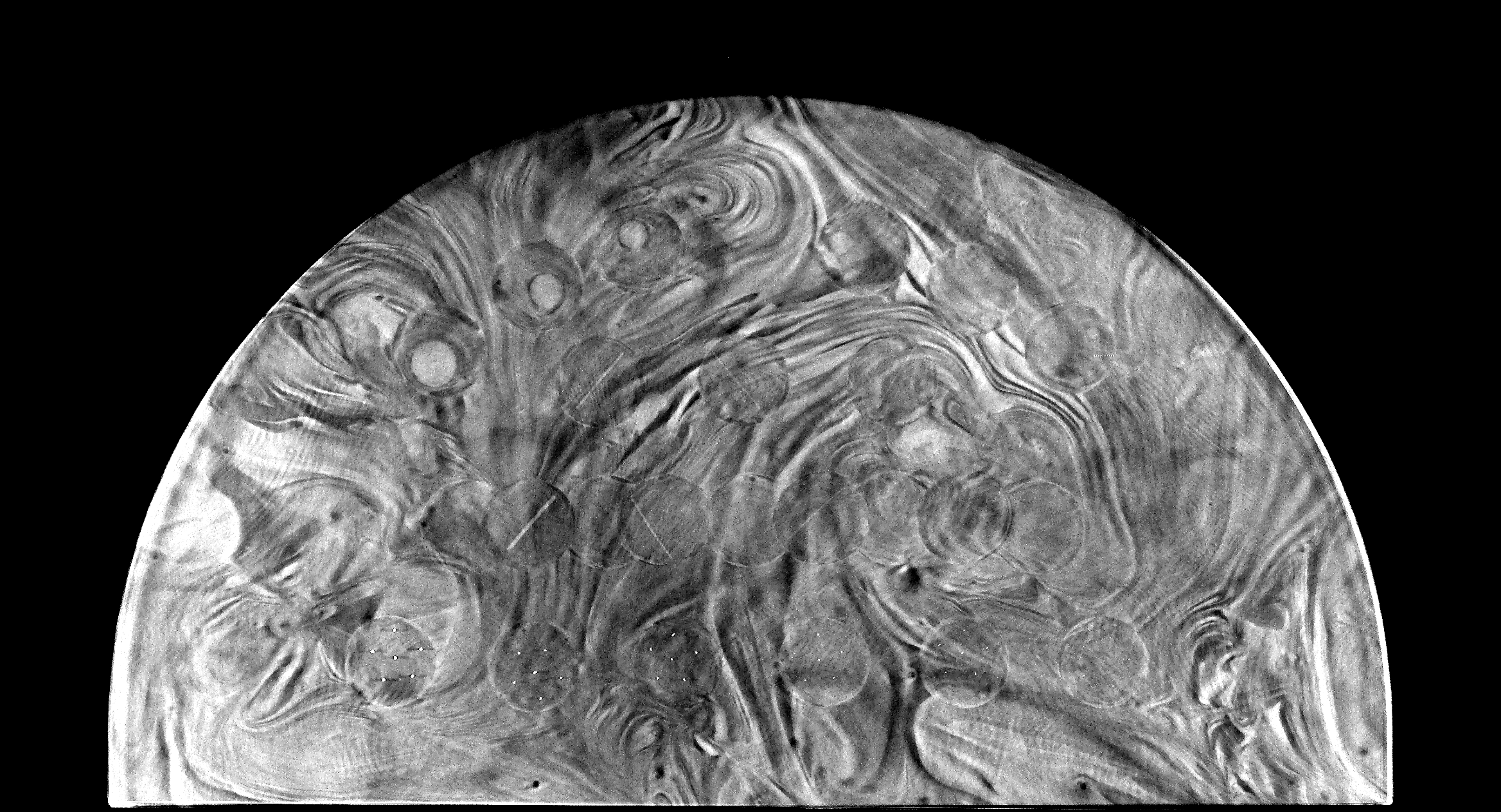} \label{fig:MicroMortadella30}}
\subfloat[30 iterations]{\includegraphics[height=0.21\textwidth, trim={27cm 26.3cm 50.5cm 13.5cm}, clip]{Sgp_RECO_29_Layer7.png} \label{fig:MassaMortadella30}}
\hfil
\caption{SGP results on BR3D phantom. (a)-(c) Reconstructions of  MC cluster number 5  obtained after 5, 15 and 30 iterations. 
(f)-(h) Reconstructions of  mass number 2  obtained after 5, 15 and 30 iterations. }
\label{fig:SgpConv_reco}
\end{figure}

\begin{figure*}[!t] 	
\centering
\subfloat[Plane Profile ]{\includegraphics[width=0.47\textwidth]{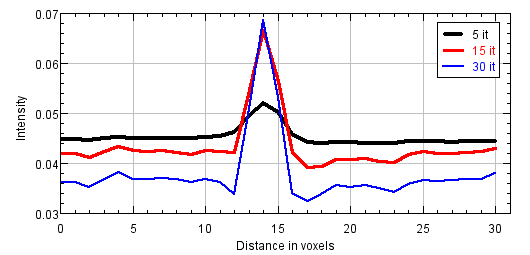} \label{fig:Mortadella_VPmicro}} 
\subfloat[Plane Profile] {\includegraphics[width=0.47\textwidth]{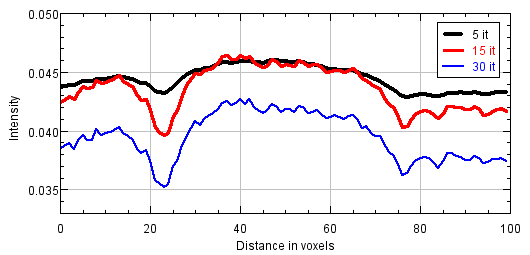} \label{fig:Mortadella_VPmassa}} 
\hfil \\
\subfloat[ASF] {\includegraphics[width=0.47\textwidth]{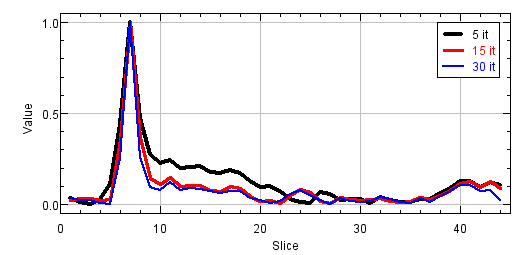} \label{fig:Mortadella_DPmicro}}
\subfloat[ASF] {\includegraphics[width=0.47\textwidth]{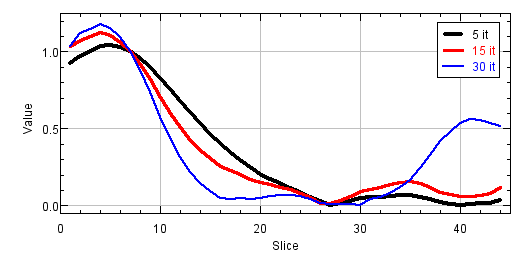} \label{fig:Mortadella_DPmassa}} 
\hfil
\caption{SGP results on BR3D phantom.  (d)-(e) Plane and Depth profile
on one microcalcification of cluster 5. 
(i)-(j) Plane Profiles and ASF profile  on the mass. In all the plots: black line corresponds to  5 iterations, red line to 15 iterations and blue line to 30 iterations.}
\label{fig:curveSgpConv}
\end{figure*}

In fact, we run the SGP solver on the BR3D phantom until the stopping condition    
 \begin{equation}
   \left |  \frac{f(x^{(k)})-f(x^{(k-1)})}{f(x^{(k)})} \right | < 10^{-6}  \quad 
     \label{eq:stop}
 \end{equation}
is satisfied. It occurs after 44 iterations. 
In Figure \ref{fig:ObjFuncMortadella}, we  plot the objective function values vs. the number of iterations: we observe that the objective function fast decreases in the first 5 iterations, whereas it exhibits a very flat trend from 10 iterations on, as it is confirmed by the red labelled values.
We have seen, in fact, that  the reconstructed images are visually almost indistinguishable after  30 iterations.\\
In Figure \ref{fig:SgpConv_reco} we exhibit the reconstructions of the 165 $\mu$m MCs of cluster 5, and the 4.7 mm mass (MS 2), obtained by the SGP algorithm after 5, 15 and 30 iterations.
In Figure \ref{fig:curveSgpConv} we report the corresponding PP and ASF plots.
Figure \ref{fig:SgpConv_reco} (a) shows that the MC of cluster 5 can be clearly visible after only 5 iterations and  the PP plots of Figure  \ref{fig:curveSgpConv} (a) confirms that the it gets more and more enhanced from the background. The ASF plot in Figure  \ref{fig:SgpConv_reco} (c) shows an improvement in the object detection along the $Z$ direction.
As visible in Figure \ref{fig:SgpConv_reco} and from the PP plot of Figure  \ref{fig:curveSgpConv} (b), MS 2 is out of focus at 5 iterations but its contours are  more and more defined when the algorithm approaches to convergence.
The previous plots confirm that the proposed model with TV regularization is more effective in recovering high contrast objects such as microcalcifications than low absorbing structures such as masses.

In Table \ref{tab:CNR} we report the values of the CNR parameter on the examined reconstructed microcalcifications and masses.
In particular, recalling the CNR definition \eqref{eq:cnrms} for the masses,  the background area is a circle with diameter of 80 voxels, whereas we have considered circles of diameter 40  and 25  voxels inside the masses 2 and 4, respectively.
When we compute the CNR value on a MC with equation \eqref{eq:cnrmc} we consider  the background as a circle of 20 voxels diameter and we compute $\mathcal{M}$ on a small circle of diameter 5 voxels containing the microcalcification.
Table \ref{tab:FWHM} shows the values of the FWHM index defined in \eqref{eq:fwhm} and computed on one of the reconstructed microcalcification in each cluster. The corresponding MCs width $w$ computed as in \eqref{eq:w} in micrometers are reported to be compared to the values of the diameters of the actual objects, shown in Table \ref{tab:SpecificheMortadella}. 
Both tables demonstrate that we can get improved and more accurate reconstructions as  the SGP approaches to the convergence: the increasing CNR indexes exhibit good denoising effects whereas the object enhancement is confirmed by the FWHM decreasing values.
We remark that the MCs of cluster 6 are not discernible from the background in only 5 iterations (the FWHM is not measurable on the sixth MC cluster), because they are 130 $\mu$m width and they should approximately fill inside only  two voxels. 
However,  they  can be well recovered after more iterations with  a good approximation of their real size.

\begin{table}[h]
\centering
\begin{tabular}{c| c c c }
& \multicolumn{3}{c}{CNR} \\
		                & 5 it. & 15 it.  & 30 it. \\
\hline 
MC cluster 3			& 24.21     & 33.34     & 38.00 \\
MC cluster 5	    	& 10.03     & 19.00     & 28.00 \\
MC cluster 6 	    	& 7.27      & 11.02     & 17.00 \\ 
\hline
MS 2		            & 0.82      & 1.07      & 1.66  \\
MS 4		            & 0.87      & 1.00      & 1.33  \\
\end{tabular}
\caption{Values of the CNR index computed after 5, 15 and 30 SGP iterations. The CNR value computed on microcalcifications is defined as in \eqref{eq:cnrmc}, whereas the CNR value computed on the masses is defined as in \eqref{eq:cnrms}.
}
\label{tab:CNR}
\end{table}
\begin{table}[h]
\centering
\begin{tabular}{c| c c c | c c c }
MC    & \multicolumn{3}{c|}{FWHM} & \multicolumn{3}{c}{ $w$ ($\mu$m)} \\
cluster	& 5 it. & 15 it.  & 30 it. & 5 it. & 15 it.  & 30 it. \\
\hline
     3	 &4.77   &3.32   &2.70   & 430   &299	&243     \\
     5	 &3.52	& 2.65  &2.32   & 317   &238    &209    \\
     6	 &-  &2.05   &1.52   & -  &185    &137     \\
\end{tabular}
\caption{ FWHM index \eqref{eq:fwhm} and  $w$ measures \eqref{eq:w} computed on the reconstructed MCs of the BR3D phantom, after 5, 15 and 30 SGP iterations. }
\label{tab:FWHM}
\end{table}


\subsection{Experiments on a human data set \label{subs:results_seno}}
\begin{figure*}[!t] 	
\centering
\subfloat[5 iterations]{\includegraphics[height=0.3\textwidth, angle =90]{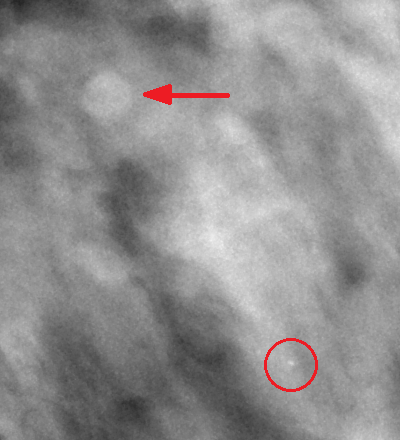} \label{fig:Seno4_fetta29}}
\subfloat[15 iterations]{\includegraphics[height=0.3\textwidth, angle =90]{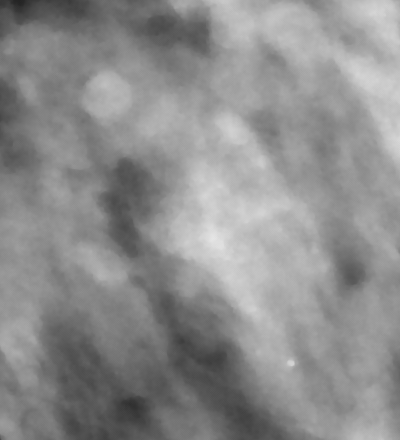} \label{fig:Seno12_fetta29}}
\subfloat[30 iterations]{\includegraphics[height=0.3\textwidth, angle =90]{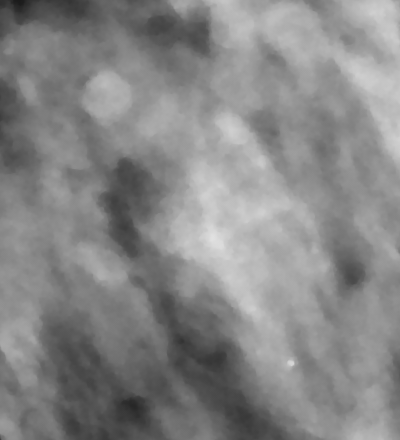} \label{fig:Seno30_fetta29}}
\hfill  \\
\subfloat[Plane Profile on the mass]{\includegraphics[width=0.47\textwidth]{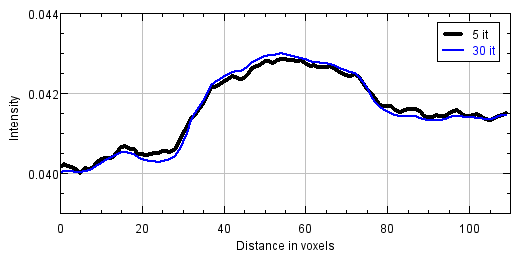} \label{fig:Seno_massa}} 
\subfloat[Plane Profile on the MC]{\includegraphics[width=0.47\textwidth]{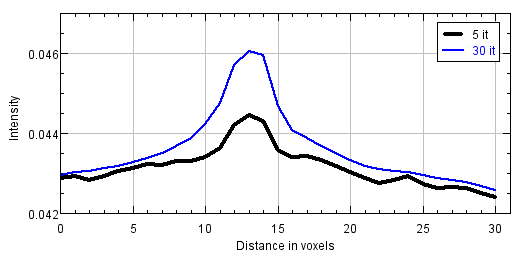} \label{fig:Seno_micro}}
\caption{Results obtained after 5, 15 and 30 SGP iterations on a human breast data set.
(a)-(c) Reconstructions of a 440 $\times$ 400 pixels region presenting both a spherical mass (pointed by the arrow) and a microcalficication (identified by the circle).
(d)-(e)  Plane profiles on the mass and on the microcalcification. In  the plots: black line corresponds to  5 iterations and blue line to 30 iterations.}
\label{fig:recoSeno}
\end{figure*}

We now illustrate the results obtained by reconstructing a real breast volume with the SGP solver at different iterative stages.
In Figure \ref{fig:recoSeno}  (a)-(c)  we report a crop of a reconstructed slice, where we can distinguish objects of interest, i.e.  a spherical mass and a small microcalcification.
The plots in Figure \ref{fig:recoSeno} (d)-(e) represent the PP calculated on the mass and the microcalcification, respectively. 
The mass is well distinguishable since the earliest reconstruction and its shape and gray level intensity do not change remarkably; however the regular blue and thin profile in Figure \ref{fig:recoSeno} (d)  points out the denoising effects of the TV function in the last iterations.
Also the microcalcification is detected in few iterations, even if a more time-consuming SGP execution enhances the contrast of the object with respect to the background.
Table \ref{tab:ValoriSeno}, reporting CNR and FWHM values computed on the objects in Figure \ref{fig:recoSeno}, gives more insight on the quality of the reconstructions.
In particular, it  confirms that the noise progressively decreases and the microcalcification gets more and more defined, from 5 to 30 iterations.

\begin{figure*}[!t] 	
\centering
\subfloat[5 iterations]{\includegraphics[height=0.30\textwidth, angle =90]{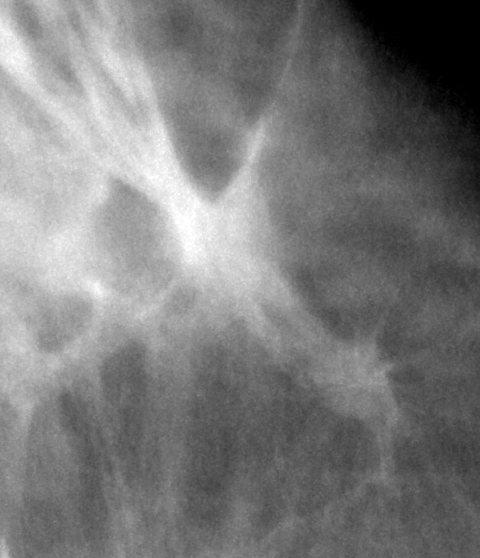} \label{fig:Seno4_fetta42}}
\subfloat[15 iterations]{\includegraphics[height=0.30\textwidth, angle =90]{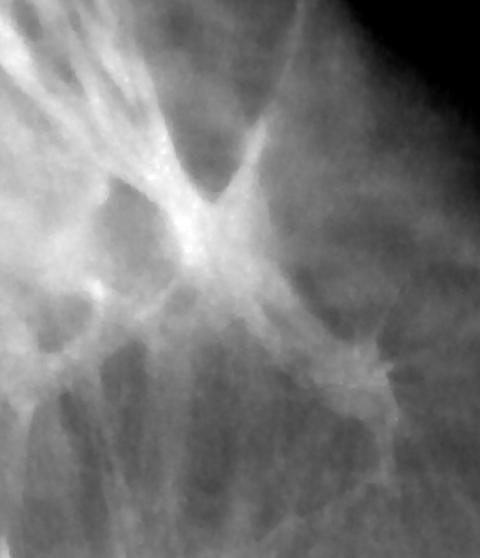} \label{fig:Seno12_fetta42}}
\subfloat[30 iterations]{\includegraphics[height=0.30\textwidth, angle =90]{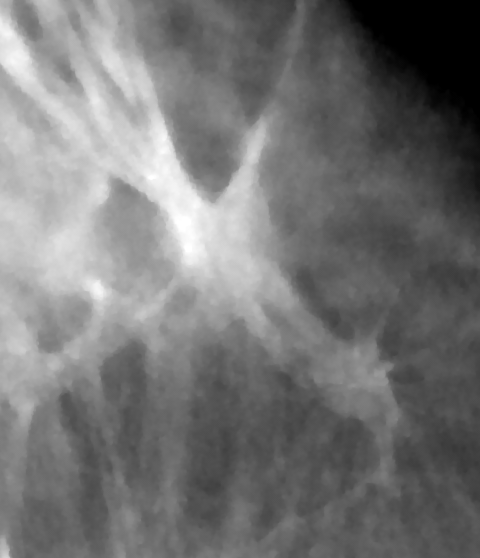} \label{fig:Seno30_fetta42}}
\hfill 
\caption{Results obtained after 5, 15 and 30 SGP iterations on a human breast data set. The reported   558 $\times$ 480 pixels crops  present two spiculated masses.
 }
\label{fig:recoSenoSpiculato}
\end{figure*}

In  Figure \ref{fig:recoSenoSpiculato}, we report the reconstruction of two  spiculated masses, which can occur in clinical cases. For such breast objects the previous measures of merits are not applicable. However we can observe that they  both are  well recognizable in the earliest reconstruction and the edges become sharper with increasing iterations.

\begin{table}[h]
\centering
\begin{tabular}{c| c c c | c c c }
    & \multicolumn{3}{c|}{CNR} & \multicolumn{3}{c}{FWHM}  \\
      & 5 it. & 15 it.  & 30 it. & 5 it. & 15 it.  & 30 it. \\
\hline 
MS   	& 0.239      & 0.381      & 0.558   & -         & -         & - \\
MC		& 8.78       & 16.59      & 16.49   & 8.57      & 7.81      & 7.29\\
\end{tabular}
\caption{Values of CNR and FWHM measures on the mass and the microcalcification  observable in Figure \ref{fig:recoSeno} (a)-(c).}
\label{tab:ValoriSeno}
\end{table}

\subsection{Experiments with a variable regularization parameter  \label{subs:paramreg}}

In all the above experiments, we have set a constant value of the regularization parameter along the iterations, to let the solvers perform at their best on the prefixed model derived by the settled $\lambda$.
In this paragraph we show the results achieved with the automatic rule \eqref{eq:UpdateLambda} for the choice of a decreasing sequence $\{ \lambda_k\}_k$ applied to the SGP solver, to reconstruct the BR3D phantom. 

\begin{figure}[!t] 	
\centering 
\includegraphics[width=0.43\textwidth]{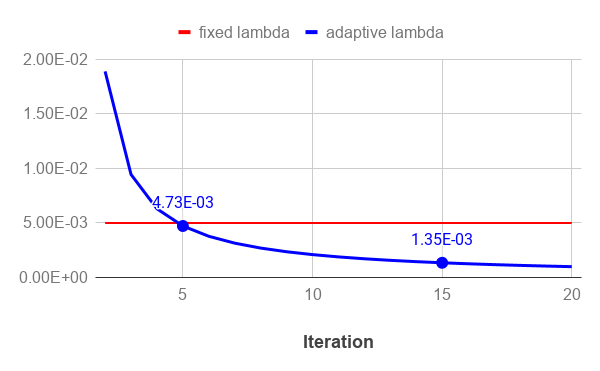} 
\caption{Sequence of decreasing $\lambda_k$ values versus the number of iterations (blue line) in the  SGP execution on the phantom test. The red straight line represents  the constant value $\lambda = 0.005$ used in SGP for the experiments presented in the previous sections. }
\label{fig:LambdaVar}
\end{figure}

\begin{figure}[!t] 	
\centering 
\subfloat[5 iterations]{\includegraphics[width=0.47\textwidth]{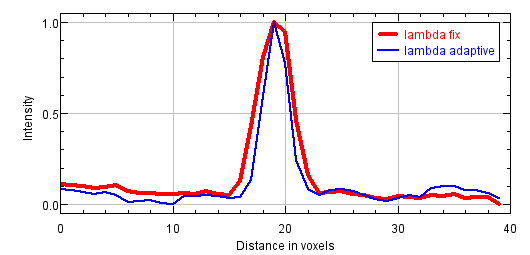} \label{fig:LambdaVar_5iter}}\\
\subfloat[15 iterations]{\includegraphics[width=0.47\textwidth]{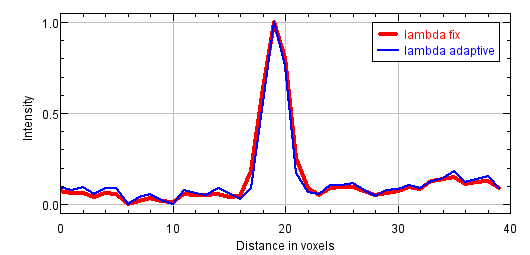} \label{fig:LambdaVar_15iter}}
\caption{Plane Profiles on one microcalcification of cluster number 3 of the phantom, obtained with SGP with different regularization parameters, in 5 and 15 iterations.
In all the plots: red line corresponds to  fix parameter, blue line the adaptive choice of $\lambda$.}
\label{fig:LambdaVariabile}
\end{figure}

Figure \ref{fig:LambdaVar} plots the sequence $\{ \lambda_k\}_k$  with the blue line, while the red line represents the constant $\lambda$ value used in the SGP implementation in the previous experiments.  
We observe that the proposed strategy  computes values greater than the  heuristically fixed one $\lambda = 0.005$ until the fifth iteration.
In Figure \ref{fig:LambdaVariabile} we compare the PP  of one microcalcification from cluster 3, reconstructed at 5 and 15 iterations using both a fixed  value  and the proposed strategy for  the regularization parameter.
We can infer from Figure \ref{fig:LambdaVariabile} (a) that the resulting larger TV weights in the first iterations  produce more accurate results. However, on advanced reconstructions the differences are negligible.
We can conclude that the proposed automatic strategy results very efficient.

\section{Conclusions \label{concl}} 


In this paper, we have presented a general optimization framework including a TV regularized formulation for DBT image reconstruction. 
We have also proposed a user-independent rule for selecting suitable values of the regularization parameter.

The results obtained with three solvers are encouraging. In early reconstructions, objects of interest of size greater than 150 $\mu m$ are visible and correctly located in the volume, whereas the object detection quality improves and the noise drastically reduces  if more iterations are allowed.
When extending the computation from 5 to 30 algorithms iterations, the increasing rate of the CNR value lies in a range $+150 \%$ to $+280 \%$.
At last, we have shown that  varying the regularization parameter  along the iterations produces better results, especially in the early stage of the algorithm execution, when compared to the use of  a fixed value heuristically chosen.\\
Since the three considered solvers produce comparable high  quality  reconstructions, we can conclude that the proposed optimization problem statement can be successfully used to detect the most interesting objects in an early diagnosis of breast tumor.

\section{Acknowledgments}
The research has been funded by the Indam GNCS grant 2020 \em{Ottimizzazione per l'apprendimento automatico e apprendimento automatico per l'ottimizzazione}.

 \bibliography{DBTRecons}

\end{document}